\documentclass[a4paper]{jpconf}

\usepackage{graphicx}

\usepackage{amsmath}
\usepackage{amsthm}
\usepackage{amssymb}
\usepackage{paralist}
\usepackage{graphics}
\usepackage{epsfig}
\usepackage[colorlinks=true]{hyperref}
\hypersetup{urlcolor=blue, citecolor=red}
\usepackage{algorithmic}
\usepackage{url}

\begin{document}
\title{Mathematical Modeling and Solution of the Tomography Problem in Domains with Reflecting Obstacles}

\author{Kamen Lozev}

\address{3161 S. Sepulveda Boulevard, Apt. 307, Los Angeles, CA 90034}

\ead{kamen@ucla.edu}

\begin{abstract}
This work develops new numerical methods for the solution of the tomography problem in domains with reflecting obstacles. We compare the solution's performance 
for Lambertian reflection, for classical tomography with ubroken rays and for specular reflection. Our numerical method using Lambertian reflection improves 
the solution's accuracy by an order of magnitude compared to classical tomography with ubroken rays and for tomography in the presence of a specularly reflecting 
obstacle the numerical method improves the solution's accuracy approximately by a factor of three times. We present efficient new algorithms for the solution's software 
implementation and analyze the solution's performance and effectiveness.  
\end{abstract}

\section{Introduction}

\bigskip

Let $f(x)$ be a continuous function in $\Omega$, where $\Omega=\Omega_0 \backslash  \overline{\Omega_1}$ and $\Omega_0$ is a 
compact convex set in $\mathbb{R}^2$ with a smooth boundary and $\Omega_1$ a convex obstacle with a smooth boundary such that
$\overline{\Omega_1} \subset \Omega_0 \subset \mathbb{R}^2$. 
Consider $\partial{\Omega_0}$ as the observation boundary of $\Omega_0$ and the points of this boundary 
as transmitters and receivers of ray signals or ray solutions of the wave equation 
\begin{equation}\label{waveequation}
u_{tt} - c^2(x)\Delta{u}=0
\end{equation}

where $c(x)>0$ is the variable speed of sound in $\Omega$  and $ u |_{\partial{\Omega_1}} = 0 $.     

The data for the tomography problem are all integrals $\int_{\gamma} f(l) dl = C_{\gamma}$ where $\gamma$ are  
rays in $\Omega$ that have both of their end points in $\partial{\Omega_0}$. One end point is a transmitter 
and the other end point of each of the rays a receiver. The classical tomography problem is to find $f(x)$ in 
$\Omega$ knowing the values of all integrals $C_{\gamma}$ where $\gamma$
are all straight line segments or unbroken rays in $\Omega$. 
In the case of a domain $\Omega_0$ without an obstacle $\Omega_1 \subset \Omega_0$, this problem is 
widely studied theoretically and numerically \cite{Na1,KS,NW,F}. When there are obstacles present, 
this problem is much less studied. Key theoretical work is done in \cite{E1,E2,E3,E4} for domains 
with one obstacle and with both broken rays, i.e. rays reflecting at the obstacle, and unbroken rays. 
A key result from \cite{E1} is that the tomography problem in the presence of one reflecting obstacle 
is well-posed. In other words, if an obstacle $\Omega_1 \subset \Omega_0$ is present then we have the well-posed problem of recovering 
$f(x)$  in $\Omega_0 \backslash \Omega_1$ from the set $C_{\gamma}$ where $\gamma$ are all broken and unbroken 
rays in the domain starting and ending at the observation boundary. This problem is called the Broken Ray Tomography Problem\cite{E1}. 

In a basic tomography setup transmitters and receivers of wave signals are placed at the domain's boundary $\partial{\Omega_0}$. 
Ray signals are generated by the transmitters and received by the receivers. Travel times for signal propagation from transmitters to 
receivers are measured and these travel times $T(A, B)$ are the values of line integrals of a function $f(x)=\frac{1}{c(x)}$ 
where $c(x)>0$ is the speed of sound at point $x \in \Omega_0 \backslash \Omega_1$. This measurement procedure gives the $C_{\gamma}$ 
data for solving the Tomography Problem by relating signal travel times to the values of line integrals of $f(x)$. 
Given sufficient data, $f(x)$ and from here the velocity $c(x)$ are computed with tomographic reconstruction algorithms \cite{Na1,KS,NW,F}. 
The tomographic algorithms presented in the next section require and compute the ray path $\gamma$ of each ray and, in order to compute these ray paths, 
we consider reflection models at $\partial{\Omega_1}$ and models of the speed of sound in $\Omega$. 

The theory of broken ray tomography from \cite{E1,E2,E3,E4} is based on a reflection model at $\partial{\Omega_1}$ that is mirror-like i.e. 
the angle of incidence is equal to the angle of reflection. In this paper, we consider a Lambertian reflection model in which incident rays are
reflected at the obstacle in all possible directions and present results that show that the broken ray tomography reconstruction 
error is smaller when we consider Lambertian reflection. 

We consider a mathematical model of the speed of sound 
\begin{equation}\label{close_to_constant_speed_of_sound}
c(x)=c_{o}+\epsilon(x)
\end{equation} 
where $x \in \Omega$ that models the speed of sound as a continuous function close to a known constant $c_{o}$. 
It is shown in \cite{RO} that for sufficiently small $\epsilon(x)$ waves propagate along the known geodesics of $c_{o}(x)$ 
when 
\begin{equation}\label{continuous_speed_of_sound}
c(x)=c_{o}(x)+\epsilon(x) 
\end{equation}
where $c(x)$, $c_{o}(x)$ and $\epsilon(x)$ are continuous functions in $\Omega$. The acoustic geodesics for constant speed of sound 
$c_o(x)=c_{o}$ are straight lines when there is no obstacle. Therefore, for the model \ref{close_to_constant_speed_of_sound}, we consider two cases. 
In the first case, $\gamma = \tilde{\gamma_1}$ is an unbroken ray composed of a straight line segment 
$\tilde{\gamma_1}$. In the second case of a broken ray, $\gamma=\tilde{\gamma_1} \bigcup \tilde{\gamma_2}$ is the union of two straight 
line segments that intersect at a reflection point at the obstacle. 

For unbroken rays the travel time or time of flight is
$$ T(A,B)=\int_{\gamma}\frac{ds}{c_o + \epsilon{(x(s))}} = \int_{\tilde{\gamma_1}}\frac{ds}{c_o + \epsilon{(x(s))}} $$ and this 
model leads to the classical tomography problem with $f(x)=\frac{1}{c_o+\epsilon(x)}$. 

\bigskip

For broken rays and a known obstacle, we know that the acoustic wave $u(x)$ propagates along the known straight line segments 
$\tilde{\gamma_1}$ and $\tilde{\gamma_2}$. $\gamma_1$ and $\gamma_2$ are known because in addition to the time of flight, our data
measurement procedure gives the end points of the ray $\gamma$ and its initial velocity, which in turn imply the reflection point of 
$\gamma$ at the known obstacle $\Omega_1$. 
Then $$ T(A,B)=\int_{\gamma}\frac{ds}{c_o + \epsilon{(x(s))}}=\int_{\tilde{\gamma_1}}\frac{ds}{c_o + \epsilon{(x(s))}} + \int_{\tilde{\gamma_2}}\frac{ds}{c_o+\epsilon{(x(s))}} $$
and, when data of this type is added to the set of measurements for the time of flight for unbroken rays, this gives the set 
$C_{\gamma}$ for the Broken Ray Tomography Problem with $f(x)=\frac{1}{c_o+\epsilon(x)}$. 
In other words, the data set $C_{\gamma}$ for the Broken Ray Tomography Problem contains the travel times of all broken and unbroken 
rays in the domain that start and end at the observation boundary.
\bigskip

\section{Numerical Solution of the Broken Ray Tomography Problem for Lambertian Reflection}

The first algorithm for the Broken Ray Tomography Problem is presented in \cite{L1} 
for a known obstacle $\Omega_1$, specular reflection and the model \ref{close_to_constant_speed_of_sound} for the speed of sound. 
This work extends the first numerical solution of the Broken Ray Tomography Problem and develops an 
algorithm for finding the velocity structure of $\Omega$ for Lambertian reflectance at $\partial{\Omega_1}$. The broken ray tomography 
problem for an obstacle with Lambertian reflectance is solved with the following algorithm that constructs the finite set of broken and 
unbroken rays  and computes the associated ray travel times by numerical integration or from ray travel time data.

\begin{algorithmic}
\REQUIRE Domain $\Omega_0$
\REQUIRE Obstacle $\Omega_1 \subset \Omega_0$
\REQUIRE Finite set of receiver points $R=(x_r, y_r)$ on $\partial{\Omega_0}$
\REQUIRE Finite set of transmitter points $T=(x_t, y_t)$  on $\partial{\Omega_0}$
\REQUIRE Finite set of obstacle boundary points $H=(x_h, y_h)$ on $\partial{\Omega_1}$
\REQUIRE Number of broken rays $n_b$. 
\REQUIRE Number of unbroken rays $n_u$.

\COMMENT{Algorithm for reconstructing $f(x,y)$ in $\Omega$ in the presence of obstacle $\Omega_0$}

\STATE{initialize an empty list L that will contain all broken and unbroken rays}

\FOR{$p = 1 \to n_b$}
  \STATE{generate a random broken ray r from input data R, T, H, $\Omega_0$ and $\Omega_1$}
  \COMMENT{a broken ray is a unique triple $((x_r, y_r), (x_t, y_t), (x_h, y_h))$ that does not intersect the obstacle except at the reflection point $(x_h, y_h)$}
  \STATE{add r to L}
\ENDFOR

\FOR{$p = 1 \to n_u$}
  \STATE{generate a random unbroken ray r from input data R, T, H, $\Omega_0$ and $\Omega_1$}
  \COMMENT{an unbroken ray is a unique pair $((x_r, y_r), (x_t, y_t))$ that does not intersect the obstacle}
  \STATE{add r to L}
\ENDFOR

\STATE{randomize L as a preprocessing step before starting the Kaczmarz method}

\FORALL{rays r in L}
\STATE{Compute travel time $p_r$ for ray r. In numerical simulations, this is obtained via nimerical integration in $\Omega_0$  along r. Store $p_r$ in list P.}  
\ENDFOR

\STATE{Reconstruct $f(x,y)$ by the Kaczmarz method \cite{K} for the linear system associated with rays r in L and corresponding travel times $p_r$ in P.}  

\end{algorithmic}

As input to the above algorithm give $n_b$ and $n_u$ to be equal or approximately equal to the maximum number of broken and unbroken rays for the finite 
input sets R, T and H, or alternatively, during ray generation, generate new rays until all rays with end points in the input sets are generated. 
In addition, again in order to approximate the requirements of the theory of broken ray tomography for inclusion of all broken and unbroken rays, provide as input 
as many transmitter and receiver points on the observation boundary as possible. Ideally, the transmitter and receiver points are uniformly 
distributed on the boundary and sufficiently close to each other so that each set approximates the set of all points on the observation boundary.
The above algorithm uses Lambertian reflection by considering broken rays as random triples $((x_r, y_r), (x_t, y_t), (x_h, y_h))$ 
and uses all rays in the domain. 

\section{Experimental Results}

In order to show the effectiveness of the numerical solution of the broken ray tomography problem with Lambertian reflection, a Java implementation of the above
algorithm compares the reconstructed values of $f(x,y)$ in $\Omega_0$ with the known values of $f(x,y)$ for the same test function $f(x,y)=K\sqrt{(x-x_0)^2+(y-y_0)^2}$ 
, where $(x_0,y_0)$ is the center of the computation grid, and test environment as in \cite{L1}. The following table compares the reconstruction error 
for classical tomography without reflection and broken ray tomography with Lambertian reflection.

\begin{table}[h]
\begin{tabular}{|c|cc|cc|}
	\hline
Experiment & ART Error  & ART Iterations & BRTL Error  & BRTL Iterations  \\
  \hline

1 & 1.269592e-004 & 54293 & 1.272610e-005 & 32635 \\ 
2 & 1.300814e-004 & 54204 & 5.272083e-006 & 40557 \\ 
3 & 1.478464e-004 & 45923 & 7.530267e-006 & 37251 \\ 
4 & 1.642026e-004 & 42480 & 1.110909e-005 & 33104 \\ 
5 & 1.927561e-004 & 23037 & 1.163345e-005 & 33424 \\ 
6 & 1.985439e-004 & 33705 & 1.896151e-005 & 28522 \\ 
7 & 8.991362e-005 & 88190 & 1.401511e-005 & 30061 \\ 
8 & 1.641201e-004 & 40717 & 1.955676e-005 & 29687 \\ 
9 & 1.089771e-004 & 64932 & 7.102756e-006 & 40126 \\ 
10 & 1.934379e-004 & 32899 & 2.149935e-005  & 28723 \\ 

  \hline

Average & 1.516838e-004 & 48038 & 1.294065e-005 & 33409 \\

  \hline
\end{tabular}
\caption{Error and number of iterations for broken ray tomography with Lambertian reflection(BRTL) at the boundary of the reflecting obstacle 
for a fixed number of 126050 rays. The average error for BRTL is 1.294065e-005 and the average number 
of iterations for finding a solution is 33409, and are shown in the right two columns of the table. 
The results for classical tomography with the ART method with the same number of rays are shown in the left two columns of the table. 
The average error for ART is 1.516838e-004 and the average number of iterations is 48038.
\label{BRTLambertianData1}}
\end{table}

The results show that the reconstruction error of the new numerical solution of the broken ray tomography problem using Lambertian reflection is 
approximately three times smaller compared to an average reconstruction error using specular reflection of 3.874848e-005 
by the algorithm from \cite{L1} for the same test function, domain and obstacle and an order of magnitude smaller than the reconstruction 
error for classical tomography in the presence of a reflecting obstacle. The new method also appears to be faster compared to classical tomography with ART.
The order of the rays in both the ART and BRTL tests in this work is randomized before applying the Kaczmarz method therefore the speed of convergence
difference is due to the use of reflection. We will report further results on the speed of convergence.

\section{Acknowledgements}

I would like to thank Professor Gregory Eskin for his continuous guidance. 

\section*{References}

\bibliography{BRT}

\end{document}